\documentclass[a4paper,twoside]{amsart}
\usepackage{amsmath,amssymb,amsfonts,amsthm}

\newtheorem{The}{Theorem}[section]
\newtheorem{Pro}[The]{Proposition}

\theoremstyle{definition}
\newtheorem{Def}[The]{Definition}
\newtheorem{Rem}[The]{Remark}

\title{{Second order tangent bundles of infinite dimensional manifolds}}
\author{C.T.J. Dodson}
  \address{Department of  Mathematics \\
    UMIST \\
    Manchester, M60 1QD, United Kingdom}
  \email{dodson@umist.ac.uk}
\author{G.N. Galanis}
  \address{Naval Academy of Greece \\
    Xatzikyriakion, Piraeus 185 39   \\
    Greece}
  \email{ggalanis@cc.uoa.gr}
\date{}
\begin{document}
\maketitle

\begin{abstract}
The second order tangent bundle $T^{2}M$ of a smooth manifold $M$ consists of
the equivalent classes of curves on $M$ that agree up to their acceleration.
It is known~\cite{Dod} that in the case of a finite $n$-dimensional manifold $M$, $T^{2}M$
becomes a vector bundle over $M$ if and only if $M$ is endowed with a linear connection.
Here we extend this result to $M$ modeled on an arbitrarily chosen Banach space and more
generally to those Fr\'{e}chet manifolds which can be obtained as projective
limits of Banach manifolds. The result may have application in the study of
infinite-dimensional dynamical systems.
\end{abstract}

\section*{Introduction}

The notion of the second order tangent bundle $T^{2}M$ of a smooth manifold $%
M$ as the equivalent classes of curves on $M$ that agree up
to their acceleration, seems to be a natural generalization of the classical
notion of tangent bundle $TM$ of $M$. However, the definition of a vector
bundle structure on $T^{2}M$ not only is not as evident as in the case of
tangent bundles but, in fact, is not always possible.

Dodson and Radivoiovici~\cite{Dod} proved that in the
case of a finite $n$-dimensional manifold $M$, a vector bundle structure on $%
T^{2}M$ can be well defined if and only if $M$ is endowed with a linear
connection. More precisely, $T^{2}M$ becomes then and only then a vector bundle over $M$
with structural group the general linear group $GL(2n;\mathbb{R})$ and,
therefore, a $3n$-dimensional manifold.

In this paper, we extend the aforementioned results to a wide class of
infinite dimensional manifolds. First we consider a manifold $M$ modeled
on an arbitrarily chosen Banach space $\mathbb{E}$. Using the Vilms~\cite{Vi}
point of view for connections on infinite dimensional vector bundles and a
new formalism, we generalize Dodson and Radivoiovici's main theorem by
proving that $T^{2}M$ can be thought of as a Banach vector bundle over $M$
with structural group $GL(\mathbb{E}\times \mathbb{E})$ if and only if $M$
admits a linear connection.

Taking one step further, we study also the case of Fr\'{e}chet (non-Banach)
modeled manifolds. In this framework things proved much more complicated
since there are intrinsic difficulties with Fr\'{e}chet spaces. For example,
pathological general linear groups, which do not even admit reasonable topological group
structures, put in question even the way of defining vector bundles.
However, by restricting ourselves to those Fr\'{e}chet manifolds which can
be obtained as projective limits of Banach manifolds (see e.g. \cite{Gal1}),
it is possible to endow $T^{2}M$ with a vector bundle structure over $M$
with structural group a new topological (and in a generalized sense Lie)
group which replaces the pathological general linear group of the fiber
type. This construction is equivalent with the existence on $M$ of a
specific type of linear connection characterized by a generalized set of
Christoffel symbols.

The result should in principle be of interest in the study of
infinite-dimensional dynamical systems, since important geometrical and physical
properties are normally associated with curvature properties of trajectories
in the system state space and this curvature is controlled by the second order
tangent structure. The new result provides conditions on when the space of
accelerations is simplified by the existence of a connection.

\section{Preliminaries}

In this Section we summarize all the necessary preliminary material that we
need for a self contained presentation of our paper.

Let $M$ be a $C^{\infty }-$manifold modeled on a Banach space $\mathbb{E}$
and $\{(U_{\alpha },\psi _{\alpha })\}_{\alpha \in I}$ a corresponding
atlas. The latter gives rise to an atlas $\{(\pi _{M}^{-1}(U_{\alpha }),\Psi
_{\alpha })\}_{\alpha \in I}$ of the tangent bundle $TM$ of $M$ with%
\begin{equation*}
\Psi _{\alpha }:\pi _{M}^{-1}(U_{\alpha })\longrightarrow \psi _{\alpha
}(U_{\alpha })\times \mathbb{E}:[c,x]\longmapsto (\psi _{\alpha }(x),(\psi
_{\alpha }\circ c)^{\prime }(0)),
\end{equation*}%
where $[c,x]$ stands for the equivalence class of a smooth curve $c$ of $M$
with $c(0)=x$ and $(\psi _{\alpha }\circ c)^{\prime }(0)=[d(\psi _{\alpha
}\circ c)(0)](1)$. The corresponding trivializing system of $T(TM)$ is
denoted by $\{(\pi _{TM}^{-1}(\pi _{M}^{-1}(U_{\alpha })),\widetilde{\Psi }%
_{\alpha })\}_{\alpha \in I}$.

Adopting the formalism of Vilms~\cite{Vi}, a connection on $M$ is a
vector bundle morphism:%
\begin{equation*}
D:T(TM)\longrightarrow TM
\end{equation*}%
with the additional property that the mappings $\omega _{\alpha }:\psi
_{\alpha }(U_{\alpha })\times \mathbb{E}\rightarrow \mathcal{L}(\mathbb{E},%
\mathbb{E)}$ defined by the local forms of D:%
\begin{equation*}
D_{\alpha }:\psi _{\alpha }(U_{\alpha })\times \mathbb{E}\times \mathbb{E}%
\times \mathbb{E}\rightarrow \psi _{\alpha }(U_{\alpha })\times \mathbb{E}
\end{equation*}%
with $D_{\alpha }:=\Psi _{\alpha }\circ D\circ (\widetilde{\Psi }_{\alpha
})^{-1},$ $\alpha \in I,$ via the relation
\begin{equation*}
D_{\alpha }(y,u,v,w)=(y,w+\omega _{\alpha }(y,u)\cdot v),
\end{equation*}%
are smooth. Furthermore, $D$ is a linear connection on $M$ if and only if $%
\{\omega _{\alpha }\}_{\alpha \in I}$ are linear with respect to the second
variable.

Such a connection $D$ is fully characterized by the family of Christoffel
symbols $\{\Gamma _{\alpha }\}_{\alpha \in I}$ , which are smooth mappings%
\begin{equation*}
\Gamma _{\alpha }:\psi _{\alpha }(U_{\alpha })\longrightarrow \mathcal{L}(%
\mathbb{E},\mathcal{L}(\mathbb{E},\mathbb{E}))
\end{equation*}%
defined by $\Gamma _{\alpha }(y)[u]=\omega _{\alpha }(y,u)$, $(y,u)\in \psi
_{\alpha }(U_{\alpha })\times \mathbb{E}$.

The requirement that a connection be well defined on the common areas of charts
of $M$, leads the Christoffel symbols satisfying the following compatibility
condition:

(1)$\ \ \ \left\{
\begin{array}{c}
\Gamma _{\alpha }(\sigma _{\alpha \beta }(y))(d\sigma _{\alpha \beta
}(y)(u))[d(\sigma _{\alpha \beta }(y))(v)]+(d^{2}\sigma _{\alpha \beta
}(y)(v))(u)= \\
=d\sigma _{\alpha \beta }(y)((\Gamma _{\beta }(y)(u))(v)),%
\end{array}%
\right\} $ \ \ \ \

for all $(y,u,v)\in \psi _{\alpha }(U_{\alpha })\times \mathbb{E}\times
\mathbb{E}$, and $d$, $d^{2}$ stand for the first and the second differential
respectively. For further details and the relevant proofs we refer to \cite{Vi}.

In the sequel we give some hints for a class of Fr\'{e}chet manifolds that
we will employ in the last Section of this note. Let $\{M^{i};\varphi
^{ji}\}_{i,j\in \mathbb{N}}$ be a projective system of Banach manifolds
modeled on the Banach spaces $\{\mathbb{E}^{i}\}$ respectively. If we assume
that

\textbf{(i)} the models form also a projective limit $\mathbb{F}=\varprojlim
\mathbb{E}^{i}$,

\textbf{(ii) }for each $x=(x^{i})\in M$ there exists a projective system of
local charts $\{(U^{i},\psi ^{i})\}_{i\in \mathbb{N}}$ such that $x^{i}\in
U^{i}$ and the corresponding limit $\varprojlim U^{i}$ is open in $M$,

then the projective limit $M=\varprojlim M^{i}$ can be endowed with a Fr\'{e}%
chet manifold structure modeled on $\mathbb{F}$ via the charts $%
\{(\varprojlim U^{i},\varprojlim \psi ^{i})\}$. Moreover, the tangent bundle
$TM$ of $M$ is also endowed with a Fr\'{e}chet manifold structure of the
same type modeled on $\mathbb{F}\times \mathbb{F}$. The local structure now
is defined by the projective limits of the differentials of $\{\psi ^{i}\}$
and $TM$ turns out to be an isomorph of $\varprojlim TM^{i}$. Here we adopt the
definition of Leslie~\cite{LE1},~\cite{LE2} for
the differentiability of mappings between Fr\'{e}chet spaces. However, the
differentiability proposed by Kriegl and P. Michor~\cite{Mich} is also
suited to our study.

\section{Tangent bundles of order two for infinite dimensional Banach
bundles}

Let $M$ be a smooth manifold modeled on the infinite dimensional Banach
space $\mathbb{E}$ and $\{(U_{\alpha },\psi _{\alpha })\}_{\alpha \in I}$ a
corresponding atlas. For each $x\in M$ we define the following equivalence
relation on $C_{x}=\{f:(-\varepsilon ,\varepsilon )\rightarrow M$ $|$ $f$
smooth and $f(0)=x$, $\epsilon >0\}$:%
\begin{equation}
f\approx _{x}g\Leftrightarrow f^{^{\prime }}(0)=g^{\prime }(0)\text{ and }%
f^{\prime \prime }(0)=g^{\prime \prime }(0),  \tag{2}
\end{equation}%
where by $f^{^{\prime }}$ and $f^{^{\prime \prime }}$ we denote the first
and the second, respectively, derivatives of $f$:%
\begin{eqnarray*}
f^{\prime } &:&(-\varepsilon ,\varepsilon )\rightarrow TM:t\longmapsto
\lbrack df(t)](1) \\
f^{\prime \prime } &:&(-\varepsilon ,\varepsilon )\rightarrow
T(TM):t\longmapsto \lbrack df^{\prime }(t)](1).
\end{eqnarray*}

\begin{Def}
We define the \textit{tangent space of order two} of $M$ at the point $x$ to be the
quotient $T_{x}^{2}M=C_{x}/\approx _{x}$ and the \textit{tangent bundle of order
two} of $M$ the union of all tangent spaces of order
2: $T^{2}M:=\underset{x\in M}{\cup }T_{x}^{2}M$.
\end{Def}

It is worth noting here that $T_{x}^{2}M$ can be always thought of as a
topological vector space isomorphic to $\mathbb{E}\times \mathbb{E}$ via the
bijection
\begin{equation*}
T_{x}^{2}M\overset{\simeq }{\longleftrightarrow }\mathbb{E}\times \mathbb{E:}%
[f,x]_{2}\longmapsto ((\psi _{\alpha }\circ f)^{\prime }(0),(\psi _{\alpha
}\circ f)^{\prime \prime }(0)),
\end{equation*}%
where $[f,x]_{2}$ stands for the equivalence class of $f$ with respect to
$\approx _{x}$. However, this structure depends on the choice of the chart $%
(U_{\alpha },\psi _{\alpha })$, hence a definition of a vector bundle
structure on $T^{2}M$ cannot be achieved by the use of the aforementioned
bijections. The most convenient way to overcome this obstacle is to assume
that the manifold $M$ is endowed with an additional structure: a linear
connection.

\begin{The}\label{T2Mvb}
If we assume that a linear connection $D$ is defined on the manifold $M$,
then $T^{2}M$ becomes a Banach vector bundle with structural group the
general linear group $GL(\mathbb{E}\times \mathbb{E)}$.
\end{The}

\begin{proof}
Let $\pi _{2}:T^{2}M\rightarrow M$ be the natural projection of $T^{2}M$ to $%
M$ with $\pi _{2}([f,x]_{2})=x$ and $\{\Gamma _{\alpha }:\psi _{\alpha
}(U_{\alpha })\longrightarrow \mathcal{L}(\mathbb{E},\mathcal{L}(\mathbb{E},%
\mathbb{E}))\}_{a\in I}$ the Christoffel symbols of the connection $D$ with
respect to the covering $\{(U_{a},\psi _{a})\}_{a\in I}$ of $M$. Then, for
each $\alpha \in I$, we define the mapping $\Phi _{\alpha }:\pi
_{2}^{-1}(U_{\alpha })\longrightarrow U_{\alpha }\times \mathbb{E}\times
\mathbb{E}$ with
\begin{equation*}
\Phi _{\alpha }([f,x]_{2})=(x,(\psi _{\alpha }\circ f)^{\prime }(0),(\psi
_{\alpha }\circ f)^{\prime \prime }(0)+\Gamma _{\alpha }(\psi _{\alpha
}(x))((\psi _{\alpha }\circ f)^{\prime }(0))[(\psi _{\alpha }\circ
f)^{\prime }(0)]).
\end{equation*}%
These are obviously well defined and injective mappings. They are also
surjective since any element $(x,u,v)\in U_{\alpha }\times \mathbb{E}\times
\mathbb{E}$ can be obtained through $\Phi _{\alpha }$ as the image of the
equivalence class of the smooth curve%
\begin{equation*}
f:\mathbb{R}\rightarrow \mathbb{E}:t\mapsto \psi _{\alpha }(x)+tu+\frac{t^{2}%
}{2}(v-\Gamma _{\alpha }(\psi _{\alpha }(x))(u)[u]),
\end{equation*}%
appropriately restricted in order to take values in $\psi _{\alpha
}(U_{\alpha })$. On the other hand, the projection of each $\Phi _{\alpha }$
to the first factor coincides with the natural projection $\pi _{2}:$ $%
pr_{1}\circ \Phi _{\alpha }=\pi _{2}$. Therefore, the trivializations $%
\{(U_{\alpha },\Phi _{\alpha })\}_{a\in I}$ define a fiber bundle structure
on $T^{2}M$ and we need now to focus on the behaviour of the
mappings $\Phi _{\alpha }$ on areas of $M$ that are covered by common
domains of different charts. Indeed, if $(U_{\alpha },\psi _{\alpha }),$ $%
(U_{\beta },\psi _{\beta })$ are two such charts, let $(\pi
_{2}^{-1}(U_{\alpha }),\Phi _{\alpha })$, $(\pi _{2}^{-1}(U_{\beta }),\Phi
_{\beta })$ be the corresponding trivializations of $T^{2}M$. Taking into
account the compatibility condition (1) is satisfied by the Christoffel
symbols $\{\Gamma _{\alpha }\}$ we see that:%
\begin{equation*}
(\Phi _{\alpha }\circ \Phi _{\beta }^{-1})(x,u,v)=\Phi _{\alpha }([f,x]_{2}),
\end{equation*}%
where $(\psi _{\beta }\circ f)^{\prime }(0)=u$ and $(\psi _{\beta }\circ
f)^{\prime \prime }(0)+\Gamma _{\beta }(\psi _{\beta }(x))(u)[u]=v$. As a
result,%
\begin{equation*}
(\Phi _{\alpha }\circ \Phi _{\beta }^{-1})(x,u,v)=
\end{equation*}%
\begin{equation*}
((\psi _{\alpha }\circ \psi _{\beta }^{-1})(\psi _{\beta }(x)),d(\psi
_{\alpha }\circ \psi _{\beta }^{-1}\circ \psi _{\beta }\circ
f)(0)(1),d^{2}(\psi _{\alpha }\circ \psi _{\beta }^{-1}\circ \psi _{\beta
}\circ f)(0)(1,1)+
\end{equation*}%
\begin{equation*}
\Gamma _{\alpha }((\psi _{\alpha }\circ \psi _{\beta }^{-1})(\psi _{\beta
}(x)))(d(\psi _{\alpha }\circ \psi _{\beta }^{-1}\circ \psi _{\beta }\circ
f)(0)(1))[d(\psi _{\alpha }\circ \psi _{\beta }^{-1}\circ \psi _{\beta
}\circ f)(0)(1)]=
\end{equation*}%
\begin{equation*}
(\sigma _{\alpha \beta }(\psi _{\beta }(x)),d\sigma _{\alpha \beta }(\psi
_{\beta }(x))(u),d\sigma _{\alpha \beta }(\psi _{\beta }(x))(d^{2}(\psi
_{\beta }\circ f)(0)(1,1)
\end{equation*}%
\begin{equation*}
+d^{2}\sigma _{\alpha \beta }(\psi _{\beta }(x))(u)[u]+\Gamma _{\alpha
}(\sigma _{\alpha \beta }(\psi _{\beta }(x)))(d\sigma _{\alpha \beta }(\psi
_{\beta }(x))(u))[d\sigma _{\alpha \beta }(\psi _{\beta }(x))(u)]=
\end{equation*}%
\begin{equation*}
(\sigma _{\alpha \beta }(\psi _{\beta }(x)),d\sigma _{\alpha \beta }(\psi
_{\beta }(x))(u),d\sigma _{\alpha \beta }(\psi _{\beta }(x))(d^{2}(\varphi
_{\beta }\circ f)(0)(1,1)+\Gamma _{\beta }(\psi _{\beta }(x))(u)[u])=
\end{equation*}%
\begin{equation*}
=(\sigma _{\alpha \beta }(\psi _{\beta }(x)),d\sigma _{\alpha \beta }(\psi
_{\beta }(x))(u),d\sigma _{\alpha \beta }(\psi _{\beta }(x))(v)),
\end{equation*}%
where by $\sigma _{\alpha \beta }$ we denote the diffeomorphisms $\psi
_{\alpha }\circ \psi _{\beta }^{-1}$ of $\mathbb{E}$. Therefore, the
restrictions to the fibers
\begin{equation*}
\Phi _{\alpha ,x}\circ \Phi _{\beta ,x}^{-1}:\mathbb{E}\times \mathbb{%
E\rightarrow E}\times \mathbb{E:}(u,v)\longmapsto (\Phi _{\alpha }\circ \Phi
_{\beta }^{-1})|_{\pi _{2}^{-1}(x)}(u,v)
\end{equation*}%
are linear isomorphisms and the mappings:%
\begin{equation*}
T_{\alpha \beta }:U_{\alpha }\cap U_{\beta }\rightarrow \mathcal{L}(\mathbb{E%
}\times \mathbb{E},\mathbb{E}\times \mathbb{E}):x\longmapsto \Phi _{\alpha
,x}\circ \Phi _{\beta ,x}^{-1}
\end{equation*}%
are smooth since $T_{\alpha \beta }=(d\sigma _{\alpha \beta }\circ \psi
_{\beta })\times (d\sigma _{\alpha \beta }\circ \psi _{\beta })$ holds for
each $\alpha ,\beta \in I$.

As a result, $T^{2}M$ is a vector bundle over $M$ with fibers of type $%
\mathbb{E}\times \mathbb{E}$ and structural group $GL(\mathbb{E}\times
\mathbb{E})$. Moreover, $T^{2}M$ is isomorphic to $TM\times TM$ since both
bundles are characterized by the same cocycle $\{(d\sigma _{\alpha \beta
}\circ \psi _{\beta })\times (d\sigma _{\alpha \beta }\circ \psi _{\beta
})\}_{\alpha ,\beta \in I}$ of transition functions. \bigskip
\end{proof}

\begin{Rem}
Note that in the case of finite dimensional manifolds
the vector bundle structure obtained in the previous Theorem coincides with
that defined by Dodson and Radivoiovici, since the corresponding
transition functions are identical (see \cite{Dod}; Corollary 2), although
we based ours on a different--totally coordinate free--approach.
\end{Rem}

We conclude this Section by proving that the converse also of Theorem~\ref{T2Mvb} holds:

\begin{The}\label{convT2Mvb}
Let $M$ be a smooth manifold modeled on the Banach space $\mathbb{E}$. If
the second order tangent bundle $T^{2}M$ of $M$ admits a vector bundle
structure, with fibers of type $\mathbb{E}\times \mathbb{E}$, isomorphic to
the product of vector bundles $TM\times TM$, then a linear connection can be
defined on $M$.
\end{The}

\begin{proof}
Let $\{(\pi _{2}^{-1}(U_{\alpha }),\Phi _{\alpha })\}_{\alpha \epsilon I}$
be a trivializing cover of $T^{2}M$ which, according to the hypothesis,
restricted to the fibers $\pi _{2}^{-1}(x)\simeq \pi _{M}^{-1}(x)\times \pi
_{M}^{-1}(x)$ will have the form: $\Phi _{\alpha ,x}=\Phi _{\alpha
,x}^{1}\times \Phi _{\alpha ,x}^{2}$, where $\Phi _{\alpha ,x}^{1}$ and $%
\Phi _{\alpha ,x}^{2}$ will be linear isomorphisms from $\pi _{M}^{-1}(x)$
to $\mathbb{E}$. Then, we may construct a chart $(U,\psi _{\alpha })$ of $M$
such that $d_{x}\psi _{\alpha }(f^{^{\prime }}(0))=\Phi _{\alpha
}^{1}([f,x]_{2})$. Indeed, if $(U,\psi )$ is an arbitrarily chosen chart of $%
M$ with $U\subseteq U_{\alpha }$, we may define $\psi _{\alpha }$ as the
composition of $\psi $ with $\Phi _{\alpha ,x}^{1}\circ (d_{x}\psi )^{-1}$.
Based on these charts we define the Christoffel symbols of the desired
connection as follows:%
\begin{equation*}
\Gamma _{\alpha }(y)(u,u)=\Phi _{\alpha }^{2}([f,x]_{2})-(\psi _{\alpha
}\circ f)^{\prime \prime }(0),
\end{equation*}%
where $f$ is the curve of $M$ that generates the vector $u$ with respect to
the chart $\psi _{\alpha }$. The remaining values of $\Gamma _{\alpha }(y)$
on elements of the form $(u,v)$ with $u\neq v$ are automatically defined if
we demand $\Gamma _{\alpha }(y)$ to be bilinear. These mappings satisfy
the necessary compatibility condition (1) since the trivializations $\{(\pi
_{2}^{-1}(U_{\alpha }),\Phi _{\alpha })\}_{\alpha \epsilon I}$ agree, via
the transition functions of $T^{2}M$, on all common areas of their domains,
and, thus, give rise to a linear connection on $M$.
\end{proof}

\section{The Fr\`{e}chet case}

As we saw in the previous Section, the definition of a vector bundle
structure on the tangent bundle of order two is always possible for Banach
modeled manifolds endowed with a linear connection. However, if we take one
step further by considering a manifold $M$ modeled on a Fr\`{e}chet
(non-Banach) space $\mathbb{F}$, then things prove to be much more complicated
due to intrinsic difficulties with these types of topological vector spaces.

More precisely, the pathological structure of the general linear groups $GL(%
\mathbb{F})$, $GL(\mathbb{F}\times \mathbb{F})$ which do not even admit
non-trivial topological group structures raises the question of whether
{\em any} possible vector bundle structure can be defined on $T^{2}M$.
On the other hand, the
fact that the space of continuous linear mappings between Fr\`{e}chet spaces
does not remain in the same category of topological vector spaces, as well
as the lack of a general solvability theory of differential equations on $%
\mathbb{F}$, turns the study of connections of the manifold $M$\ into a very
complicated issue.

In this Section, by employing a new methodology--which has already been
proven successful for classical tangent and frame bundles (see \cite{Gal2}, %
\cite{VG})--we develop a vector bundle structure for the
second order tangent bundles of a certain class of Fr\`{e}chet manifolds:
those which can be obtained as projective limits of Banach manifolds.

To this end, let $M$ be a smooth manifold modeled on the Fr\`{e}chet space $%
\mathbb{F}.$ Taking into account that the latter {\em always} can be
realized as a projective limit of Banach spaces $\{\mathbb{E}%
^{i};\rho ^{ji}\}_{i,j\in \mathbb{N}}$ (i.e. $\mathbb{F\cong }\varprojlim
\mathbb{E}^{i}$) we assume that the manifold itself is obtained as the limit
of a projective system of Banach modeled manifolds $\{M^{i};\varphi
^{ji}\}_{i,j\in \mathbb{N}}$ in the sense described in the Preliminaries.
We obtain:

\begin{Pro}
The second order tangent bundles $\{T^{2}M^{i}\}_{i\in \mathbb{N}}$ form
also a projective system with limit (set-theoretically) isomorphic to $%
T^{2}M $.
\end{Pro}

\begin{proof}
For any pair of indices $(i,j)$ with $j\geq i$, we define the mapping:%
\begin{equation*}
g^{ji}:T^{2}M^{j}\rightarrow T^{2}M^{i}\text{ }\mathbb{:}\text{ }%
[f,x]_{2}^{j}\longmapsto \lbrack \phi ^{ji}\circ f,\phi ^{ji}(x)]_{2}^{i},
\end{equation*}%
where the brackets $[\cdot ,\cdot ]_{2}^{j}$, $[\cdot ,\cdot ]_{2}^{i}$
denote the classes of the equivalence relation (2) on $M^{j}$, $M^{i}$
respectively. We easily check that $g^{ji}$ is always well-defined, since
two equivalent curves $f$, $g$ on $M^{j}$ will give
\begin{equation*}
d^{(n)}\phi ^{ji}(f^{(n)}(0))=(\phi ^{ji}\circ f)^{(n)}(0)=(\phi ^{ji}\circ
g)^{(n)}(0)=d^{(n)}\phi ^{ji}(g^{(n)}(0)),\text{ \ \ }n=0,1,2,
\end{equation*}%
where $d^{(1)}\phi ^{ji}:TM^{j}\rightarrow TM^{i}$ stands for the first and $%
d^{(2)}\phi ^{ji}:T(TM^{j})\rightarrow T(TM^{i})$ for the second
differential of $\phi ^{ji}$.

\qquad On the other hand, the relations $g^{ik}\circ g^{ji}=g^{jk}$ $(j\geq
i\geq k)$, readily obtained from the corresponding ones for $\{\varphi
^{ji}\}_{i,j\in \mathbb{N}}$, ensures that $\{T^{2}M^{i};g^{ji}\}_{i,j\in
\mathbb{N}}$ is a projective system. Based now on the canonical projections $%
\phi ^{i}:M\rightarrow M_{i}$ of $M$, we define
\begin{equation*}
F^{i}:T^{2}M\rightarrow T^{2}M^{i}:[f,x]_{2}\rightarrow \lbrack \phi
^{i}\circ f,\phi ^{i}(x)]_{2}^{i}\ \ \ \ \ \ \ (i\in \mathbb{N}).
\end{equation*}%
Since $g^{ji}\circ F^{j}=F^{i}$ holds for any $j\geq i$, we obtain the
mapping
\begin{equation*}
F=\varprojlim F^{i}:T^{2}M\rightarrow \varprojlim
(T^{2}M^{i}):[f,x]_{2}\rightarrow ([\phi ^{i}\circ f,\phi
^{i}(x)]_{2}^{i})_{i\in \mathbb{N}}.
\end{equation*}%
This is an injection because $F([f,x])=F([g,x])$ gives%
\begin{equation*}
d^{(n)}\phi ^{i}(f^{(n)}(0))=(\phi ^{i}\circ f)^{(n)}(0)=(\phi ^{i}\circ
g)^{(n)}(0)=d^{(n)}\phi ^{i}(g^{(n)}(0)),\ \ n=0,1,2,
\end{equation*}%
and, therefore, $f^{(n)}(0)=g^{(n)}(0)$ $\ (n=0,1,2)$ since $TM\equiv
\varprojlim TM^{i}$ and $T(TM)\equiv \varprojlim T(TM^{i})$ with
corresponding canonical projections $\{d\phi ^{i}\}_{i\in \mathbb{N}}$ and $%
\{d^{(2)}\phi ^{i}\}_{i\in \mathbb{N}}$ respectively.

On the other hand, $F$ is also surjective since for any element $\emph{a}%
=([f^{i},x^{i}]_{2}^{i})_{i\in \mathbb{N}}\in \varprojlim (T^{2}M^{i})$ we
see that:%
\begin{equation}
\lbrack \phi ^{ji}\circ f^{j},\phi
^{ji}(x^{j})]_{2}^{i}=[f^{i},x^{i}]_{2}^{i},\text{ for }j\geq i,  \tag{3}
\end{equation}%
thus $x=(x^{i})\in M=\varprojlim M^{i}$. Moreover, if $(U=\varprojlim
U^{i},\psi =\varprojlim \psi ^{i})$ is a projective limit chart of $M$
through $x$ and $(\pi _{M}^{-1}(U)=\varprojlim \pi _{M^{i}}^{-1}(U^{i}),\Psi
=T\psi =\varprojlim T\psi ^{i})$, $(\pi _{TM}^{-1}(\pi
_{M}^{-1}(U))=\varprojlim \pi _{TM^{i}}^{-1}(\pi _{M^{i}}^{-1}(U^{i})),%
\widetilde{\Psi }=T(T\psi )=\varprojlim T(T\psi ^{i}))$ the corresponding
charts of $TM$, $T(TM)$ respectively,we obtain:%
\begin{equation*}
((\psi ^{i}\circ \phi ^{ji}\circ f^{j})(0),T\psi ^{i}((\phi ^{ji}\circ
f^{j})^{\prime }(0)))=((\psi ^{i}\circ f^{i})(0),T\psi ^{i}((f^{i})^{\prime
}(0)))\Rightarrow
\end{equation*}%
\begin{equation*}
\Rightarrow (\rho ^{ji}((\psi ^{j}\circ f^{j})(0)),T\psi ^{i}(T\phi
^{ji}((f^{j})^{\prime }(0))))=((\psi ^{i}\circ f^{i})(0),T\psi
^{i}((f^{i})^{\prime }(0))).
\end{equation*}%
As a result, the elements $u=((\psi ^{i}\circ f^{i})(0))_{i\in \mathbb{N}}$,
$v=((\psi ^{i}\circ f^{i})^{\prime }(0))_{i\in \mathbb{N}}$ belong to $%
\mathbb{F\cong }\varprojlim \mathbb{E}^{i}$. Similarly, relations (3) ensure
that\ $(\phi ^{ji}\circ f^{j})^{\prime \prime }(0)=$ $(f^{i})^{\prime \prime
}(0)$ which via the charts of $T(TM)$ defined above give $T(T\psi
^{i})((\phi ^{ji}\circ f^{j})^{\prime \prime }(0))=T(T\psi
^{i})((f^{i})^{\prime \prime }(0))$ or, equivalently, $\rho ^{ji}((\psi
^{j}\circ f^{j})^{\prime \prime }(0))=(\psi ^{i}\circ f^{i})^{\prime \prime
}(0)$, for $j\geq i$. Therefore, $w=((\psi ^{i}\circ f^{i})^{\prime \prime
}(0))_{i\in \mathbb{N}}$ belongs also to $\mathbb{F\cong }\varprojlim
\mathbb{E}^{i}$. Considering now the curve $h$ of $\mathbb{F}$ with $%
h(t)=u+t\cdot v+\frac{t^{2}}{2}\cdot w$, as well as the corresponding one $f$
of $M$ with respect to the chart $(U=\varprojlim U^{i},\psi =\varprojlim
\psi ^{i})$, we may check that
\begin{equation*}
(\phi ^{i}\circ f)(0)=\phi ^{i}(x)=x^{i}=f^{i}(0),
\end{equation*}%
\begin{eqnarray*}
(\phi ^{i}\circ f)^{\prime }(0) &=&(\psi _{i}^{-1}\circ \rho _{i}\circ
h)^{\prime }(0)=T\psi _{i}^{-1}((\rho _{i}\circ h)^{\prime }(0))= \\
&=&T\psi _{i}^{-1}(\rho _{i}(v))=T\psi _{i}^{-1}((\psi _{i}\circ
f^{i})^{\prime }(0))= \\
&=&(f^{i})^{\prime }(0),
\end{eqnarray*}%
\begin{eqnarray*}
(\phi ^{i}\circ f)^{\prime \prime }(0) &=&(\psi _{i}^{-1}\circ \rho
_{i}\circ h)^{\prime \prime }(0)=T(T\psi _{i}^{-1})((\rho _{i}\circ
h)^{\prime \prime }(0))= \\
&=&T(T\psi _{i}^{-1})(\rho _{i}(w))=T(T\psi _{i}^{-1})((\psi _{i}\circ
f^{i})^{\prime \prime }(0))= \\
&=&(f^{i})^{\prime \prime }(0),
\end{eqnarray*}%
for all indices $i,j$ with $j\geq i$. As a result, the curves $\phi
^{i}\circ f$, $f^{i}$ are equivalent on $M^{i}$ and $%
F([f,x]_{2})=([f^{i},x^{i}]_{2}^{i})_{i\in \mathbb{N}}=\emph{a}$.

\qquad By this means, we ensure that the mapping $F$ is the desired
isomorphism which turns $T^{2}M$, $\varprojlim (T^{2}M^{i})$ to isomorphic
sets.
\end{proof}

Based on the last result, next we define a vector bundle
structure on $T^{2}M$ by means of a certain type of linear connection on $M.$
The problems concerning the structural group of this bundle (discussed
earlier)\ are overcame by the replacement of the pathological $GL(\mathbb{F}%
\times \mathbb{F})$ by the new topological (and in a generalized sense
smooth Lie) group:%
\begin{equation*}
\mathcal{H}^{0}(\mathbb{F\times F}):=\{(l^{i})_{i\in \mathbb{N}}\in {%
\prod_{i=1}^{\infty }}GL(\mathbb{E}^{i}\mathbb{\times E}^{i}):\,\varprojlim
l^{i}\,\text{\ exists}\}.
\end{equation*}%
To be more specific, $\mathcal{H}^{0}(\mathbb{F\times F})$ is a topological
group being isomorphic to the projective limit of the Banach-Lie groups
\begin{equation*}
\mathcal{H}_{i}^{0}(\mathbb{F\times F}):=\{(l^{1},l^{2},...,l^{i})_{i\in
\mathbb{N}}\in {\prod_{k=1}^{i}}GL(\mathbb{E}^{k}\mathbb{\times E}%
^{k}):\,\rho _{jk}\circ l^{j}=l^{k}\circ \rho _{jk}\,\text{\ }(k\leq j\leq
i)\}.
\end{equation*}%
On the other hand, it can be considered as a generalized Lie group via its
embedding in the topological vector space $\mathcal{L}(\mathbb{F\times F})$.
Using these notations we obtain:

\begin{The}
If a Fr\`{e}chet manifold $M=\varprojlim M^{i}$ is endowed with a linear
connection $D$ that can be also realized as a projective limit of
connections $D=\varprojlim D^{i}$, then $T^{2}M$ is a Fr\'{e}chet vector
bundle over $M$ with structural group $\mathcal{H}^{0}(\mathbb{F\times F}).$
\end{The}

\begin{proof}
Following the terminology established in Section 1, we consider $%
\{(U_{\alpha }=\varprojlim U_{\alpha }^{i},\psi _{\alpha }=\varprojlim \psi
_{\alpha }^{i})\}_{\alpha \in I}$ an atlas of $M$. Each linear connection $%
D^{i}$ $(i\in \mathbb{N)}$, which is naturally associated to a family of
Christoffel symbols $\{\Gamma _{\alpha }^{i}:\psi _{\alpha }^{i}(U_{\alpha
}^{i})\rightarrow \mathcal{L}(\mathbb{E}^{i},\mathcal{L}(\mathbb{E}^{i},%
\mathbb{E}^{i}))\}_{\alpha \in I}$, ensures that $T^{2}M^{i}$ is a vector
bundle over $M^{i}$ with fibers of type $\mathbb{E}^{i}$. This structure, as
already presented in Theorem~\ref{T2Mvb}, is defined by the trivializations:%
\begin{equation*}
\Phi _{\alpha }^{i}:(\pi _{2}^{i})^{-1}(U_{\alpha }^{i})\longrightarrow
U_{\alpha }^{i}\times \mathbb{E}^{i}\times \mathbb{E}^{i},
\end{equation*}%
with
\begin{equation*}
\Phi _{\alpha }^{i}([f,x]_{2}^{i})=(x,(\psi _{\alpha }^{i}\circ f)^{\prime
}(0),(\psi _{\alpha }^{i}\circ f)^{\prime \prime }(0)+\Gamma _{\alpha
}^{i}(\psi _{\alpha }^{i}(x))((\psi _{\alpha }^{i}\circ f)^{\prime
}(0),(\psi _{\alpha }^{i}\circ f)^{\prime }(0)));\ \alpha \in I.
\end{equation*}

Taking into account that the families of mappings $\{g^{ji}\}_{i,j\in
\mathbb{N}}$, $\{\varphi ^{ji}\}_{i,j\in \mathbb{N}}$, $\{\rho
^{ji}\}_{i,j\in \mathbb{N}}$ are connecting morphisms of the projective
systems $T^{2}M=\varprojlim (T^{2}M^{i})$, $M=\varprojlim M^{i}$, $\mathbb{E=%
}\varprojlim \mathbb{E}^{i}$ respectively, we check that the projections $%
\{\pi _{2}^{i}:T^{2}M^{i}\rightarrow M^{i}\}_{i\in \mathbb{N}}$ satisfy
\begin{equation*}
\varphi ^{ji}\circ \pi _{2}^{j}=\pi _{2}^{i}\circ g^{ji}\text{ \ }(j\geq i%
\mathbb{)}
\end{equation*}%
and the trivializations $\{\Phi _{\alpha }^{i}\}_{i\in \mathbb{N}}$%
\begin{equation*}
(\varphi ^{ji}\times \rho ^{ji}\times \rho ^{ji})\circ \Phi _{\alpha
}^{j}=\Phi _{\alpha }^{i}\circ g^{ji}\text{ \ }(j\geq i\mathbb{).}
\end{equation*}%
As a result,%
\begin{equation*}
\pi _{2}=\varprojlim \pi _{2}^{i}:T^{2}M\longrightarrow M
\end{equation*}%
exists and is a surjective mapping,
\begin{equation*}
\Phi _{\alpha }=\varprojlim \Phi _{\alpha }^{i}:\pi _{2}^{-1}(U_{\alpha
})\longrightarrow U_{\alpha }\times \mathbb{F}\times \mathbb{F}\text{ \ }%
(\alpha \in I)
\end{equation*}%
is smooth, as a projective limit of smooth mappings, and its projection to
the first factor coincides with $\pi _{2}$.

\qquad On the other hand, the restrictions of $\Phi _{\alpha }$ to any fiber
$\pi _{2}^{-1}(x)$ is a bijection since $\Phi _{\alpha ,x}:=pr_{2}\circ \Phi
_{\alpha }|_{\pi _{2}^{-1}(x)}=\varprojlim (pr_{2}\circ \Phi _{\alpha
}^{i}|_{(\pi _{2}^{i})^{-1}(x)})$.

The crucial part of our construction, however, concerns the corresponding
transition functions $\{T_{\alpha \beta }=\Phi _{\alpha ,x}\circ \Phi
_{\beta ,x}^{-1}\}_{\alpha ,\beta \in I}$. These can be considered as taking
values in the generalized Lie group $\mathcal{H}^{0}(\mathbb{F}\times
\mathbb{F)}$, since $T_{\alpha \beta }=\epsilon \circ T_{\alpha \beta
}^{\ast }$, where $\{T_{\alpha \beta }^{\ast }\}_{\alpha ,\beta \in I}\ $ are
the smooth mappings
\begin{equation*}
T_{\alpha \beta }^{\ast }:U_{\alpha }\cap U_{\beta }\rightarrow \mathcal{H}%
^{0}(\mathbb{F}\times \mathbb{F}):x\longmapsto (pr_{2}\circ \Phi _{\alpha
}^{i}|_{(\pi _{2}^{i})^{-1}(x)})_{i\in \mathbb{N}}
\end{equation*}%
and $\epsilon $ is the natural inclusion
\begin{equation*}
\epsilon :\mathcal{H}^{0}(\mathbb{F}\times \mathbb{F})\rightarrow \mathcal{L}%
(\mathbb{F\times F}):(l^{i})_{i\in \mathbb{N}}\longmapsto \varprojlim l^{i}.
\end{equation*}%
Summarizing, we have proved that $T^{2}M$ admits a vector bundle structure
over $M$ with fibers of type $\mathbb{F}\times \mathbb{F}$ and structural
group $\mathcal{H}^{0}(\mathbb{F}\times \mathbb{F})$. Moreover, this bundle
is isomorphic to $TM\times TM$ since we may check that they have identical
transition functions:
\begin{equation*}
T_{\alpha \beta }(x)=\Phi _{\alpha ,x}\circ \Phi _{\beta ,x}^{-1}=(d(\psi
_{a}\circ \psi _{\beta }^{-1})\circ \psi _{\beta })(x)\times (d(\psi
_{a}\circ \psi _{\beta }^{-1})\circ \psi _{\beta })(x)
\end{equation*}
\end{proof}

We conclude this paper by proving that also the converse of the previous Theorem
is true.

\begin{The}
If $T^{2}M$ is an $\mathcal{H}^{0}(\mathbb{F}\times \mathbb{F})-$Fr\`{e}chet
vector bundle over $M$ isomorphic to $TM\times TM$, then $M$ admits a linear
connection which can be realized as a projective limit of connections.
\end{The}

\begin{proof}
By the hypothesis, the vector bundle structure on $T^{2}M$ would be defined
by a family of trivializations $\{\Phi _{\alpha }:\pi _{2}^{-1}(U_{\alpha
})\longrightarrow U_{\alpha }\times \mathbb{F}\times \mathbb{F}\}_{\alpha
\in I}$ which will be realized as projective limits of corresponding
trivializations $\Phi _{\alpha }^{i}:(\pi _{2}^{i})^{-1}(U_{\alpha
}^{i})\longrightarrow U_{\alpha }^{i}\times \mathbb{E}^{i}\times \mathbb{E}%
^{i}$ of $T^{2}M^{i}$ $(i\in \mathbb{N})$ so that the transition functions $%
\{T_{\alpha \beta }\}_{\alpha ,\beta \in I}$ of $T^{2}M$ take their values
in $\mathcal{H}^{0}(\mathbb{F}\times \mathbb{F})$. As a result, every
factor-bundle $T^{2}M^{i}$ will be a vector bundle isomorphic to $%
TM^{i}\times TM^{i}$ and, according to Theorem 2.4, a linear connection $%
D^{i}$ can be defined on $M^{i}$ with Christoffel symbols satisfying
\begin{equation*}
\Gamma _{\alpha }^{i}(y)([f,x]^{i},[f,x]^{i})=(pr_{3}\circ \Phi _{\alpha
}^{i})([f,x]_{2}^{i})-(\psi _{\alpha }^{i}\circ f)^{\prime \prime }(0).
\end{equation*}%
We may check then that $\varprojlim (\Gamma _{\alpha }^{i}(y^{i})(u^{i}))$
exists for any $(y^{i})$, $(u^{i})\in \mathbb{F=}\varprojlim \mathbb{E}^{i}$%
, thus the connections $\{D^{i}\}_{i\in \mathbb{N}}$ form a projective
system with corresponding limit the desired linear connection $D=\varprojlim
D^{i}$ on $M$.
\end{proof}

\end{document}